\input amstex
 \input amsppt.sty
 \documentstyle{amsppt}
 \loadmsbm
 \NoBlackBoxes
 \leftheadtext {(S. Liriano)}
 \rightheadtext{( Dimension and deviation)}
  \topmatter
 \title Krull Dimension and Deviation in certain
 Parafree Groups
 \endtitle
 \author Sal Liriano
 \endauthor
 \affil
 SAL21458\@yahoo.com
 \endaffil
  \dedicatory
  In memory of
 Marcelo Llarull (1958-2005)
  \enddedicatory
 \medskip
 \abstract
 Hanna Neumann asked whether it was
 possible for two non-isomorphic residually nilpotent finitely
 generated ({\it fg})
 groups, one of them free,
 to share the lower central sequence. G. Baumslag answered the
 question in the affirmative and thus gave rise to parafree groups.
 A group $G$ is termed parafree of rank
 $n$ if it is residually nilpotent and shares the same lower central
 sequence with a free group of rank $n$. The deviation of a {\it
 fg} parafree group $G$ of rank $n$ is the difference $\mu(G)-n$, where $\mu(G)$
  is the minimum possible number of generators of $G$.
 \medskip

 Let $G$ be {\it fg;} then $Hom(G,SL(2,\Bbb C))$ inherits
 the structure of an algebraic variety, denoted by $R(G)$, which is
 an invariant of {\it
 fg}
 presentations of $G$. If $G$ is an $n$
 generated parafree group, then the deviation of
  $G$ is 0 iff $Dim(R(G))=3n$. It is known that for $n\ge 2 $ there exist
  infinitely many
  parafree groups of rank $n$ and deviation 1 with non-isomorphic
  representation varieties of dimension $3n$. In this paper it is
  shown that given integers $n \ge 2$, and $k\ge 1$, there exists
  infinitely many parafree groups of rank $n$ and deviation $k$
  with non-isomorphic representation varieties of dimension
  different from $3n$; in particular, there exist infinitely
  many parafree groups $G$ of rank $n$ with $Dim(R(G))> q$,
  where $q \ge 3n$ is an arbitrary integer.
 \endabstract
 \endtopmatter
 \document
 \remark{Structure of paper} New results in this paper are Theorem
 1, Theorem 2, Theorem 3, Theorem 4, and Theorem 5. This paper is
 broken up into 3 parts: Introduction, Section one, Section two.
 In the introduction an outline of the mentality that guided
 this investigation is given, along with the proof of some
 preliminary results including the proof of Theorem 2, and Theorem
 3. In Section One  material involving sequences of primes, and
 groups associated with such sequences is developed. The section
 ends with a proof of Theorem 4. Section Two begins with a proof
 of Theorem 1, and ends with a proof of Theorem 5.
 \endremark

 \head Introduction
 \endhead

 Let $G$ be a finitely generated group.
 Then the set of homomorphisms from $G$ into $SL_2\Bbb C$
 inherits the structure of an algebraic variety denoted by $R(G)$.
 This algebraic variety, also known as the representation variety of
 $G$, is an invariant of finitely generated presentations of
 $G$. It exports
 into the theory of finitely generated
 groups the numerous invariants of Commutative Algebra and Algebraic
 Geometry. The object of this work  is to continue the exploration
 with this invariant of a rather remarkable class of groups invented by G.
 Baumslag in the 1960's, parafree groups. For the convenience of the reader, the
 definition of a parafree group follows.

\medskip\flushpar
 Let $\gamma_n F$ stand for the $n$-th term of the Lower Central Series
 of the free group $F$. Then, a group $G$ is termed
  parafree  if:
  \roster
 \item "{1)}" G is residually nilpotent.

 \item "{2)}" There exists a free group $F$ with the property that
 $G/ \gamma_n G \cong F/ \gamma _n F $ for all $n\ge  1$.
 \endroster

\medskip\flushpar
 Now denote by $\mu(G)$ the minimal number of generators of  a
 group  $G$.  Define the {\it  rank \/ } of $G$, denoted here by $rk(G)$, to be
 $rk(G) = \mu (G/ \gamma_2 G)$. Further, define the
 {\it  deviation  \/ }
  of the group $G$, here denoted by $\delta(G)$, to be
 $\delta(G) = \mu(G)-rk(G)$. A  parafree  group $G$ is
 of rank $r$ if the free group in (2) above  is also of rank $r$.

\medskip
Perhaps two of the most fundamental invariants of an algebraic
variety are its Krull dimension and reducibility status. As it is
often the case, with these invariants can be associated other
invariants. For example, let $c$ be a
 positive integer. Now denote by $N_c(V)$
 the number of maximal irreducible components ({\it mirc}) of $V$ of dimension
 $c$. Then $N_c(V)$ is an invariant
 of the algebraic variety
 $V$. This invariant was introduced by the author in
 \cite {L3} where it was shown that  $N_4(R(\frak T))=\bold g$,
  where $\bold g$ is
 the genus of the
 torus knot corresponding to the group $\frak T$. Subsequently, in
 \cite {L4}
 the invariant $N_c(V)$ was employed in studying a class of
 parafree groups to show,
 amongst other
 surprising results, that
 given any parafree group $G$ of rank $n$, deviation 1, and with
 $Dim(R(G))=3n$, that there exist a set $P_n$ of
 parafree groups of rank $n$ and deviation 1 having the property
 that no two
 groups in $P_n$ have isomorphic representation varieties, and
 such that for any $G_j\in P_n$, it is the case that
 $N_{3n}(R(G))< N_{3n}(R(G_j))$. In fact, the set,
 $\{N_{3n}(R(G_i))\vert G_i \in P_n \}$ is infinite. This is significant
 given
 the bewildering likeness between a {\it
 fg} free group \footnote {A free group
 of rank $n$ is a parafree group of deviation 0 and rank $n$.}
 and a {\it
 fg} parafree group of the same rank with deviation other than 0.
 Indeed, a non-free {\it
 fg} parafree group of
 rank $n$ agrees with a free group of the same rank on an infinite
 number of {\it
 fg}
 torsion-free nilpotent quotients.
 The stage is now set to introduce the main results of
 this paper.

\medskip
 A direct result of W. Magnus' investigations
 is that
 a {\it
 fg} parafree group of rank $n$ is
 free iff it has deviation zero; see \cite {WM2}, \cite
 {WM}. In \cite
 {L4} it was shown that {\it an $n$-generated group
 $G$ is free, iff $Dim(R(G))=3n$}. So it follows then that:

\proclaim {Proposition 0.1}
  If $G$ is an $n$ generated parafree group
  then the deviation of
$G$ is zero iff $Dim(R(G))=3n.$ \endproclaim

 In other words, up to isomorphism a parafree group $G$ of rank $n$
 and deviation zero is determined by
an invariant of $R(G)$, namely its dimension. Admittedly, it was
Proposition 0.1 that motivated the author to pursue more closely
the possible connections between rank, deviation of a parafree
group $G$,  and the krull dimension of the algebraic variety
$R(G)$. Indeed, all example of finitely generated parafree groups
of rank $n$ considered in \cite {L4} turned out to have dimension
$R(G)=3n$, that of the free group $F_n$ with which they shared the
lower central sequence. \footnote {The {\it lower central sequence
} of a group $G$ is the sequence given by $G/\gamma_n G$, where
$\gamma_n G$ is the $n$-term of the lower central series.} In the
face of it, this seems anything but surprising given that parafree
groups share a host of properties with free groups, (see \cite
{L4}, and \cite {B1}, \cite {B2} for an account). But, since {\it
 fg} parafree groups of finite
rank can have arbitrarily large deviations, it did not seem at all
unusual to inquire whether indeed it is possible that all finitely
generated parafree groups $G$ of rank $n$ have the property that
$Dim(R(G))=3n$? The next proposition provides a bound for the
dimension of the representation variety of a parafree group of
rank $n$ and deviation $k$.

\proclaim {Proposition 0.2} Let $G$ be a parafree group of rank
$n$ and deviation $k\ge 1$, then $Dim(R(G))< 3(n+k).$
\endproclaim

\demo {Proof} If $G$ is a parafree group of rank $n$ and deviation
$k$, then it can be easily seen that $G$ is generated by $n+k$
elements. To see this, recall the formula for the deviation
$\delta(G)$ of a group introduced above: $\delta(G)=\mu(G)-rk(G)$,
where $\mu(G)$ is the smallest possible number of generators of
$G$. Inserting $k$ and $n$ as in the statement above yields
$k+n=\mu(G)$. Thus $G$, since it is not free is a quotient group
of a free group of rank $k+n$. Now, since a free group of rank
$n+k$ has a faithful representation in $SL(2,\Bbb C)$ by Sanov,
(see \cite {SN}), and since $R(F_{n+k})$ is an irreducible variety
of dimension exactly $3(n+k)$, it follows that $Dim(R(G))< Dim
(R(F_{n+k}))=3(n+k)$.
\enddemo

Next a theorem is unveiled that gives an answer of no to the
question asking whether $Dim(R(G))$ $=3n$ for any parafree group
$G$ of rank $n$. In fact, it guarantees that for each integer
$k\ge 1 $ there exist an infinite set of rank 2 parafree groups of
deviation $k$ having the property that their corresponding
representation varieties are $4+2k$-dimensional, reducible, and
pairwise not isomorphic, a fact that implies that the groups are
not isomorphic. This is an amazing illustration of the sensitivity
of some algebro-geometric invariants when unleashed in the study
of {\it
 fg} groups.

\proclaim {Theorem 1} Given an integer $k \ge 1 $ there exists an
infinite set $S_k$ consisting of parafree groups or rank 2 and
deviation $k$ having the property that for $G_i$, $G_j$ in $S_k$
with $i \ne j$ the following holds: \roster
\item $Dim (R(G_i))=4+2k$.
\item If $i\ne j$, and $c=4+2k$, then $N_c(R(G_i))\ne N_c(R(G_j)).$
\item $R(G_i)\ncong R(G_j)$.
\item Given \footnote {In particular, $m$ can be
$N_c(R(G))$ for any $G$ in $S_k$.} any integer $m$ there exists an
infinite subset $Q$ in $S_k$ such that for any $G_t \in Q$ and $c$
as immediately above $N_c(R(G_t))\ge m$.
\endroster
\endproclaim
\medskip

\medskip\flushpar The next theorem is a consequence of Theorem 1.
Because its proof  is so dependent on Theorem 1, it is given
immediately.

\proclaim {Theorem 2} Given any integer $r \ge 2$ and any integer
$k\ge 1$ there exists an infinite set $S_{r,k}$ consisting of
parafree groups of rank $r$ and deviation $k$ with the property
that for each $G_i$ and $G_j$ with $i\ne j$ in $S_{r,k}$ it is the
case that: \roster
\item  $Dim(R(G_i))=3r+2k-2.$
\item  If $c=3r+2k-2$, then $N_c(R(G_i))\ne N_c(R(G_j)).$
\item $R(G_i)\ncong R(G_j).$
\item  Given \footnote {Adjust with
$S_{r,k}$, and $c$ the statement of the footnote above.} any
integer $m$, there exists an infinite subset $Q$ in  $S_{r,k}$
such that for any $G_t \in Q$ and $c$ as directly above
$N_c(R(G_t))\ge m$.
\endroster
\endproclaim

 \demo {Proof}Let $k\ge 1$ be some
fixed arbitrary integer. If $r=2$, then replace $S_{r,k}$ by
$S_k$, as in Theorem 1. Now, suppose that $r
> 2.$ Then let
$S_{r,k}=\{F_{r-2}\ast G_t \vert G_t \in S_k \}.$ The sets
$S_{r,k}$ will satisfy all the statements of the theorem and
consist of parafree groups of rank $r$ since the free product of
two parafree groups is of rank the sum of the rank of their
factors \cite {B2}. The proof of part one is a trivial consequence
of the fact that the dimension of a product of varieties is the
sum of the dimensions of the factors, and that
$Dim(R(F_{r-2}))=3(r-2).$ Part two follows from the fact that for
$G_i$ and $G_j$ in $S_k$ with $i \ne j$ and $c=4+2k$ it is the
case that $N_c(R(G_i))\ne N_c(R(G_k))$ together with the fact that
the representation variety of a free group in an irreducible
linear algebraic group is irreducible. Three is a consequence of
the fact that the number of {\it mirc} is an invariant of the
isomorphism class of an algebraic variety. Part four is a direct
consequence of the fact that the sets $S_k$ in Theorem 1 are all
infinite and well ordered by the {\it mirc} counting function in
dimension $4+2k$, and this can be extended to $S_{r,k}$ using the
{\it mirc} counting function in dimension $3r+2k-2$. The proof is
complete.
\enddemo

\medskip
Undoubtedly, it is a well documented fact that cyclically pinched
\footnote{A group is termed cyclically pinched if it can be given
a presentation of the type $G=<X \cup Y;W=V>$, where $X$ is a
finite set of generators, $W$ is a non-trivial word in the free
group on $X$, and $Y$ is a finite set of generators, and the word
$V$ a non-trivial word in the free group on $Y$.} one-relator
 groups share a host of
properties with free groups, a fact that given the likeness
between non-free parafree groups and free groups, endows one
relator parafree groups with cyclically pinched presentations with
numerous additional free-like properties; the reader is encouraged
to consult \cite {B4}, \cite {B5}, \cite {FR}, \cite {L4} to learn
of some of them.

\medskip For $k\ge 2$ the reader will see that the groups used
in the proof of Theorem 1 turn out not to be one relator groups.
Indeed, using Proposition 0.2 it is possible to deduce that an
$n+1$ generated parafree group of rank $n$ and deviation one
always has $Dim(R(G))\le 3n+2$, and using the main theorem of
\cite {L1} the next result follows immediately.

\proclaim {Proposition 0.3} Let $G$ be a a group on $n+1$
generators and having a cyclically pinched presentation $<x_1,
\cdots, x_n, y; W=y^p >$; then $Dim(R(G)) \le 3n+1.$
\endproclaim

An immediate consequence of the above proposition is that any
parafree group obtained by adjoining a root to a non-trivial
element $W$ of the free group of rank $n \ge 2$ thus: $<x_1,
\cdots,x_n,y; W=y^p >$, which is the method described by G.
Baumslag in  \cite {B1}, can only have a representation variety of
dimension at most $3n+1$. Indeed, until now all one relator
parafree groups obtained in the described manner have resulted in
groups having representation varieties of dimension precisely
three times the rank of the free group to which the root was
adjoined. But, is this always the case? The answer is provided by
the next theorem.

\proclaim {Theorem 3} Given \footnote {In particular, $k$ can be
taken equal to $N_7(R(G))$ for any fixed {\it
 fg} parafree group $G$ with
$Dim (R(G))=7$.} integer $k \ge 1$ there exist infinitely many
parafree groups of rank 2 and deviation 1 with pairwise
non-isomorphic representation varieties of dimension precisely 7
and such that $N_7(R(G))\ge k.$
\endproclaim

\demo {Proof} The set of groups in one to one correspondence with
the set $S$ in Theorem 5 part (3) can be well ordered using the
{\it mirc} counting function $N_7(V)$ in dimension 7. Since the
set $S$ is infinite, it follows that given any integer $k$ there
exist only a finite number of integers smaller that $k$ and which
correspond to the ordering imposed by the {\it mirc} counting
function $N_7(V)$ on the set of groups corresponding to $S$.  The
result follows, since each of the groups that correspond to the
set $S$, besides being parafree of rank 1, has representation
variety of dimension 7.

\enddemo

\head Section One \endhead

\medskip
 In this section material necessary for the proof of
Theorem 1 and Theorem 5 is introduced. Because of the significant
role it will play, it is prudent to begin by recalling a result
developed and proved by the author in \cite {L4}, which in the
sequel shall be referred to as the Coarse Sieve Theorem \footnote
{The reason for the name is because under certain
algebro-geometric conditions on the representation varieties of
groups pertaining to the terms of an infinite sequence
$S_{i=1}^{\infty}(G_i,\frak N_i)$, the theorem guarantees that the
sequence $S_{i=1}^{\infty}(G_i,\frak N_i)$ has an infinite
subsequence whose corresponding groups $G_i$ have pairwise
non-isomorphic representation varieties.}, or the CS Theorem, for
short.
\medskip
\medskip

\proclaim {Coarse Sieve Theorem}
 \flushpar
 Let $S_{i=1}^{\infty}(G_i,\frak N_i)$ be an infinite
 sequence of pairs $(G_1,\frak N_1),
(G_2,\frak N_2), \cdots$ consisting of {\it fg} groups $G_i$ and a
corresponding normal subgroup $\frak N_i$ of $G_i$. Let $N_c(V)$
be the {\it mirc} counting function in dimension $c$. Suppose that
$Dim(R(G_i))=Dim(R(G_i/\frak N_i))=c$, and that the set
$\tilde{S}=\{ N_c(R(G_j/\frak N_j))\,\vert \,(G_j,\frak N_j) \in
S_{i=1} ^{\infty}(G_i,\frak N_i)\}$ contains an infinite set of
integer points. Then: \roster
\item "{i)}" $S_{i=1}^{\infty}(G_i,\frak N_i)$ has an infinite
subsequence $S^2$ with the property that given two different pairs
$(G_j,\frak N_j),\, (G_k,\frak N_k)$ in $S^2$ then
$N_c(R(G_j))\neq N_c(R(G_k)).$
\item "{ii)}" For different pairs
$(G_j,\frak N_j),\, (G_k,\frak N_k)$ in $S^2$ then $R(G_j)\ncong
R(G_k).$
\item "{iii)}" The set of groups
$S_2=\{G \vert  \text {$G$ occurs in some term of $S^2$} \}$ can
be well ordered using the function $N_c(V)$.
\endroster
\endproclaim
\medskip
Next some preparation  for the eventual deployment of the Coarse
Sieve Theorem in the proof of Theorem 5, and Theorem 1 is
introduced.
\medskip
\remark {$k$-Prime Sequences} The Sieve of Eratosthenes can be
used to find the terms of the ascending sequence of primes
starting from the prime 3  and contained in the interval $[3,m]$,
for any \footnote {For sufficiently large $m$ one of the many
estimates  of the prime counting function $\pi (m)$ may be used to
approximate the number of terms of the sequence in the interval
$[3,m]$.} arbitrary positive integer $m$. Inductively then one may
list the ascending sequence of primes starting from 3; call the
resulting sequence of primes: $1P_{i=1}^{\infty}$. Now, given
integer $k\ge 1$ use $1P_{i=1}^{\infty}$ to obtain the sequences $
2P_{i=1}^{\infty}$,..., $ kP_{i=1}^{\infty}$ specified in the
following list \footnote {Note that $1P_{i=1}^{\infty}$ may
sometimes be written as $P_{i=1}^{\infty}$ in the sequel.}:
\roster
\item "{a)}" $1P_{i=1}^{\infty}:P_1,P_2,P_3, \dots $
\item "{b)}" $2P_{i=1}^{\infty}:(P_1,P_2),(P_3,P_4),\dots$
\item "{}" .
\item "{}" .
\item "{}" .
\item "{t)}" $kP_{i=1}^{\infty}:(P_1,P_2,\dots,
P_k),(P_{k+1},\dots,P_{2k}),\dots$.
\endroster
where $P_1=3, P_2=5, P_3=7,\dots$ is the infinite sequence of
primes starting with 3, and (a) gives the sequence when $k=1$, and
(b) gives the sequence when $k=2$, and so on.

\endremark
\medskip

\medskip
Next, a procedure for associating a very particular free product
of groups with each term of any of the above $k$-prime sequences
is revealed. But first, it is important to mention that due to the
presentation of certain parafree groups of importance in the
sequel, the group $$G=<x_1,x_2;x_1^2x_2^2=1>\tag 1-1$$ shall play
a fundamental role. Fortunately, some relevant aspects of the
algebraic geometry of $R(G)$ are well understood.

\medskip
\remark {Procedure 1.1} \flushpar Let $G=<x_1,x_2;x_1^2x_2^2=1>$,
and let $k\ge 1$ be any integer. Associate with any term of the
sequence $kP_i^{\infty}$ above, say the $t$-term $kP_t$, where $t$
is any integer $\ge 1$, a unique group $GkP_t$ given by:

$$GkP_t=G\ast \Bbb Z_{P_{(t-1)k+1}}\ast \Bbb Z_{P_{(t-1)k+2}}\ast
\cdots\ast \Bbb Z_{P_{tk}} \tag 1-2$$
\endremark

\medskip
\flushpar So for example, the group associated with the third term
$P_3$ of the sequence $P_{i=1}^{\infty}$ would be $GP_3=G\ast \Bbb
Z_{P_3}$, which is nothing but the free product  $G\ast \Bbb Z_7$.
The group associated with the second term $2P_2$ belonging to the
sequence $2P_i^{\infty}$ is given by $G2P_2$, which is nothing but
$G \ast \Bbb Z_{P_3}\ast \Bbb Z_{P_4}$ which really stands for: $G
\ast \Bbb Z_7\ast \Bbb Z_{11}$.

\proclaim {Lemma 1.1} Let $G=<x_1,x_2;x_1^2x_2^2=1>$. Then
 $Dim(R(G))=4$, and $R(G)$ is reducible with $N_4(R(G))=1$.
\endproclaim

\demo {Proof} This follows from the main theorem in  \cite {L3}.
\enddemo


\proclaim {Lemma 1.2} Let $GkP_j$ be the group associated with the
$j$-th term of sequence $kP_{i=1}^{\infty}$, where $j$ is any
integer $\ge 1$. Let $c=4+2k.$ Then: \roster
\item $Dim(R(GkP_j))=c;$
\item if $k\ge 2$ then $N_c(R(GkP_j))=
\frac { (P_{(j-1)k+1}-1)\cdots (P_{jk}-1)} {2^k};$
\item if $k=1$ then $N_c(R(GkP_j))=\frac {(P_{j}-1)} {2}.$
\endroster
\endproclaim

\demo {Proof}$GkP_j$ is the free product of the group $G$ in Lemma
1.1 and a free product of cyclics involving cyclic factors that
are all of finite prime order $t \ge 3$.  The representation
variety of a cyclic group of finite order $t\ge 3$ can by
elementary methods be shown to be 2 dimensional, and when $t \ge
3$ is prime, to have $\frac {t-1} {2}$ {\it mirc} of dimension
two. See \cite {L1}, or \cite {L3}. This with the help of the next
lemma, whose proof is left as an exercise, will yield the result
sought.

\proclaim {Lemma 1.3} Let $V$ and $W$ be two algebraic varieties
with $Dim(V)=d$ with $N_d(V)=y$, and with $Dim(W)=r$ with
$N_r(W)=m$. Then $Dim(V\times W)=d+r$ and $N_{d+r}(V\times W)=ym$.
\endproclaim

\flushpar Now by Lemma 1.1, and (10), and by Lemma 1.3 it follows
that $Dim (R(GkP_j))=4+2k$ and that $N_c(R(GkP_j))= \frac {
(P_{(j-1)k+1}-1)\cdots (P_{jk}-1)} {2^k}$,\,\,  if $k\ge 2$. If
$k=1$ the result is an immediate consequence of the fore stated.
This completes the proof of Lemma 1.2.

\enddemo

As an immediate consequence of Lemma 1.2, is the following.

\proclaim {Lemma 1.4} If $k\ge 1$ is any integer, then $\lim_ {j
\to \infty} N_{4+2k}(R(GkP_j))=\infty.$
\endproclaim

\demo { Proof} See the formula for the number of
$4+2k$-dimensional maximal irreducible components in Lemma 1.2.
\enddemo

\medskip
\flushpar Yet one more procedure is necessary; it will for each
integer $k\ge 1$, and each term in the sequence
$kP_{i=1}^{\infty}$ associate a group as follows:

\remark {Procedure 1.2} For each integer $k\ge 1$ associate with
the $n$-th term, $kP_n=(P_{k(n-1)+1},$
$P_{k(n-1)+2},\cdots,P_{nk})$, of the sequence $kP_{i=1}^{\infty}$
a corresponding group $pGkP_n$ on the free generators
$\{x_1,x_2,y_1,y_2,\cdots,y_k \}$ given by the presentation:

$$<x_1,x_2,y_1,\cdots,y_k;
x_1^{2^1}x_2^{2^1}=y_1^{P_{k(n-1)+1}},
x_1^{2^2}x_2^{2^2}=y_2^{P_{k(n-1)+2}},
\cdots,x_1^{2^k}x_2^{2^k}=y_k^{P_{nk}}> \tag {1-3}$$

\flushpar {\it Example a}: if $k=1$ and $n=2$ the group associated
with the $P_2$ term of the sequence $1P_{i=1}^{\infty}$ would be
$<x_1,x_2,y_1;x_1^{2^1}x_2^{2^1}=y_1^{P_{2}}>$ which is after
evaluation the group $<x_1,x_2,y_1;x_1^{2}x_2^{2}=y_1^{5}>.$
\medskip

\flushpar {\it Example b}: if $k=2$ and $n=3$ the group associated
with the $2P_3$ term of the sequence $2P_{i=1}^{\infty}$ would be
$<x_1,x_2,y_1,y_2; x_1^{2^1}x_2^{2^1}=y_1^{P_{k(n-1)+1}},
x_1^{2^2}x_2^{2^2}=y_2^{P_{k(n-1)+2}}>$ which is after evaluation
the group $<x_1,x_2,y_1,y_2;x_1^{2}x_2^{2}=y_1^{13},
x_1^{2^2}x_2^{2^2}=y_1^{17}>.$

\endremark

\proclaim {Lemma 1.5} Given any integer $k\ge 1$ and any integer
$e\ge 1$ then $$Dim(R(pGkP_e))=Dim(R(GkP_e)), \tag {1-4}$$ where
$pGkP_e$ is the group associated with the $e$-th term of the
sequence $kP_{i=1}^{\infty}$ using Procedure 1.2, and $GkP_e$  is
the group associated with the $e$-th term of the sequence
$kP_{i=1}^{\infty}$ using Procedure 1.1.
\endproclaim

\demo {Proof} That $Dim(R(pGkP_e))=4+2k$ follows directly from
Theorem 4. That $Dim(R(GkP_e))=4+2k$ is a direct consequence of
Lemma 1.2.
\enddemo

\proclaim { Lemma 1.51} Given integer $k\ge 1$ and the group
$pGkP_j$ associated with the $j$-th term of the sequence
$kP_{i=1}^{\infty}$ then:
$$ \frac {pGkP_j} {N(x_1^2x_2^2)} \cong GkP_j \tag {1-5}$$
\flushpar where $N(x_1^2x_2^2)$ stands for the normal closure of
the word $x_1^2x_2^2$ in the group $pGkP_j$.
\endproclaim

\demo {Proof} Applying Tietze transformations yields  that:
 $$ \frac {pGkP_j} {N(x_1^2x_2^2)}=<x_1,x_2;x_1^2x_2^2=1>\ast
\Bbb Z_{P_{(j-1)k+1}}\ast \Bbb Z_{P_{(j-1)k+2}}\ast \cdots\ast
\Bbb Z_{P_{jk}}= GkP_j,$$ \flushpar the desired result.
\enddemo

At this juncture, a series of additional lemma and observations to
be employed in the proof of Theorem 4, will be introduced.

\proclaim {Lemma 1.7} Let $V$ be an algebraic variety and $W$ a
subvariety. Then $Dim(V)=\text{Maximum}\,\, \{Dim(V-W),\,
Dim(W)\}.$
\endproclaim
\demo { Proof} Obvious.
\enddemo

\medskip\flushpar
Given $a \in SL_2\Bbb C$ and $p \in \Bbb Z$, denote the set: $\{m
\vert m^p=a, m\in SL_2 \Bbb C \}$ by $\Omega (p,a)$. Further,
given a positive integer $t$ define $\Omega^t (p,a)=\Omega
(p,a)\times \cdots \times\Omega (p,a)$, \, taken $t$ times.

\proclaim {Lemma 1.71 } All matrices of $SL_2\Bbb C$ having fixed
trace $\alpha$ where $\alpha\ne \pm 2 $ form an irreducible two
dimensional variety.
\endproclaim

\demo {Proof} Any two matrices in $SL_2\Bbb C$ with the same trace
$\alpha\ne\pm 2 $ are similar and also diagonalizable. Now define
a regular map $\Gamma: SL_2\Bbb C\times n\rightarrow V$ given by
$(m,n) \rightarrow m^{-1}nm$, where $n$ is a $SL_2\Bbb C$ matrix
of trace $\alpha$ having zero in all entries that are not in its
main diagonal. The map $\Gamma$ is a regular map onto $V$, where
$V$ is the subvariety of $SL_2\Bbb C$ consisting of all matrices
of trace $\alpha$. Suppose that $V$ is reducible; then $SL_2\Bbb
C\times n$ is also reducible, a contradiction. That $V$ is two
dimensional follows from the fact that $V$ is the zero locus in
$SL_2\Bbb C$ of a non-constant polynomial, and an application of
the Krull Principal Ideal Theorem.
\enddemo

\proclaim {Corollary 1.71} $\Omega (2,-I)$ is an irreducible two
dimensional variety.
\endproclaim

\demo{Proof}
 Using the Cayley-Hamilton Theorem and the characteristic polynomial for matrices
 in $SL_2\Bbb C$ one can
deduce that $m \in \Omega (2,-I)$ iff $m$ has trace zero. Using
Lemma 1.71 above the result follows.
\enddemo

 The set of possible groups in the next theorem includes in it the
 groups $pGkP_j$, for any integer $k\ge 1$ and any integer $j\ge
1$. However, their description does not make use of the prime
sequences introduced earlier; for that reason they have been named
differently.

 \proclaim
{Theorem 4} Let $G_k= <x_1,x_2,y_1,\cdots,y_k;x_1^{2^1}x_2^{2^1}=
y_1^{p_1},\cdots,x_1^{2^k}x_2^{2^k}=y_k^{p_k}>$, where the
$p_1,\cdots, p_k$ are distinct primes all $\ge 3$, and $k$ is an
integer $\ge 1$. Then \roster
\item The groups $G_k$ are parafree of rank 2 deviation $k$;
\item $Dim(R(G_k))= 4+2k$;
\item $R(G_k)$ is reducible.
\endroster

\endproclaim

\demo{Proof}

\flushpar Proof 1: That the groups $G_k$ are parafree of rank 2
follows directly from \cite {B1}. That they are of deviation $k$
is a consequence of the Grushko-Neumann Theorem since they map
onto a free product of cyclics involving $k+2$ non-trivial
factors.

\medskip

\flushpar Proof 2 and 3:  Let $G_k=
<x_1,x_2,y_1,\cdots,y_k;x_1^{2^1}x_2^{2^1}=
y_1^{p_1},\cdots,x_1^{2^k}x_2^{2^k}=y_k^{p_k}>$, where the
$p_1,\cdots, p_k$ are distinct primes all $> 2$, and $k$ is an
integer $\ge 1$. It will be shown that $Dim(R(G_k))= 4+2k$.
\medskip\flushpar
Let $$\rho:R(G_k)\rightarrow (SL_2\Bbb C)^2 \tag {1-6}$$ be the
projection map given by $\rho(m_1,m_2, \cdots,
m_{k+1},m_{k+2})=(m_1,m_2)$. Since the $p_i$ are odd primes this
map is onto by Lemma 1.8 in \cite {L1}.

\medskip\flushpar Now for integers $i$ such that
$1\leq i\leq k$ define $S_i\subset R(G_k)$ as follows:
$S_i=\{(m_1,m_2, \cdots, m_{k+1},m_{k+2})\vert
m_1^{2^i}m_2^{2^i}=\pm I \}$, and define $S_i^+$ as members of
$S_i$ with $m_1^{2^i}m_2^{2^i}= I$; also define $S_i^-$ as those
members of $S_i$ with $m_1^{2^i}m_2^{2^i}=-I$ in an analogous
manner. Now let

$$\bold S=\cup _{i=1} ^{k} S_i ,\tag
{1-7}$$

\flushpar and define also $\bold S^+$ and $\bold S^-$ in an
analogous fashion as before. Now define

$$\rho (\bold S_i)=\{(m_1,m_2) \vert m_1^{2^i}m_2^{2^i}=\pm I, m_1,m_2
\in SL_2\Bbb C \} .$$

\flushpar Also define $\rho ( \bold S_i^+)$ and $\rho (\bold
S_i^-)$ in an analogous manner as was done before. Define $\rho (
\bold S)=\cup _{i=1} ^{k} \rho (S_i)$

\medskip\flushpar
By Lemma 1.8   \cite {L1} the regular map $\rho$ in (1-6) has the
property that $\rho ^{-1}$ restricted to the set
$$(SL_2\Bbb C)^2 - \rho (\bold S) \tag {1-8}$$
has finite fiber with coordinates over each point in (1-8) sitting
in the set
$$R(G_k)^o=R(G_k)-\bold S. \tag {1-9}$$ Thus by Proposition
2.5 in \cite {L1} it follows that
$$Dim(\overline{(R(G_k))-\bold S})=6, \tag {1-10}$$ where the
over line indicates the Zariski closure.

\medskip
Thus to compute $Dim(R(G_k))$ all that needs to be known is
$Dim(\bold S)$, by Lemma 1.7. Note that $\Omega^2 (2,-I)$ is a
subset of $\rho (\bold S)$, by Lemma 1.9 (below), and by Lemma
1.71 this is a 4 dimensional variety. So $Dim(\rho (\bold S))\ge
4$. But Theorem 0.2 of \cite {L1}, together with Lemma .3 of \cite
{L4} guarantees that $Dim(\rho (\bold S))=4$, the desired result.
\medskip

The subvariety $\rho (S_1)$ of $\rho (\bold S)$ will play a
central role as made evident by the next two lemmas, the fist of
which is obvious.

\proclaim {Lemma 1.9 } $\Omega ^2(2,-I)\subseteq \rho (S_i^+)$ for
all $i \in \{1,2,\dots, k\}$.
\endproclaim

\proclaim {Lemma 1.91 }  The fiber of $\rho$ over $\Omega
^2(2,-I)$ is of constant dimension $2k$ at each point of $\Omega
^2(2,-I)$.
\endproclaim

\demo {Proof} Let $(m_1,m_2)\in \Omega ^2(2,-I)$ then $\rho
^{-1}(m_1,m_2)=(m_1,m_2)\times \Omega (p_1,I)\times \cdots \times
\Omega (p_k,I)$ by Lemma 1.9 . Now by Lemma 1.6 of \cite {L1 },
$Dim((m_1,m_2)\times \Omega (p_1,I)\times \cdots \times \Omega
(p_k,I))=0+2k.$ Sice the dimension of a product is the sum of the
dimensions of the factors.
\enddemo

\flushpar An immediate consequence of Lemma 1.9  and Lemma 1.91
is that
$$ Dim  (\bold S) \ge 4+2k \tag {1-11}$$

\flushpar since $\bold S$ has a subvariety that maps via a regular
map onto the 4 dimensional irreducible variety $\Omega^2(2,-I)$
and has fiber of constant dimension $2k$ over each point
$(m_1,m_2) \in \Omega ^2(2,-I).$

\medskip At this juncture it is only necessary to show that
$Dim(\bold S)\le 4+2k$. To achieve this the 4 dimensional variety
$\rho ( \bold S)$ will be decomposed into its maximal irreducible
components $C_{\alpha_1}, \dots, C_{\alpha_s}$. Clearly, since
$\rho (\bold S)$ is 4-dimensional, all components are of dimension
at most 4. Further, for each $i$ it is the case that $\rho
^{-1}(C_{\alpha_i}) $ is a closed subvariety of $R(G_k)$.  Now,
without loss of generality let $(m_1,m_2)$ be any point in say
$C_{\alpha_i}$, then the fiber over $(m_1,m_2)$ has the form
$$(m_1,m_2)\times \Omega (p_1, m_1^{2^1}m_2^{2^1})\times \cdots \times
\Omega(p_k, m_1^{2^k}m_2^{2^k}), \tag {1-12}$$

\flushpar but this fiber is always of dimension at most $2k$ by
Lemma 1.6 and Lemma 1.8 of \cite {L1}. It follows then that $Dim
(\bold S) \leq 4+2k$ since the dimension of  $\rho (\bold S)$, the
base, is 4. Now, by Lemma 1.7 it follows that $Dim (R(G_k))= \text
{Maximum}\,\{6,\, 4+2k \}=4+2k$, since $k \ge 1$. The reducibility
of $R(G_k)$ follows since there is always an irreducible component
of dimension 6 and at least one irreducible component in $\bold S$
of dimension $4+2k$. This completes the proof of the theorem. The
existence of the six dimensional component is a consequence of the
following:

\proclaim {Proposition 1.1} Let $V$ be an algebraic variety with
$Dim(V)=n$, and containing an  open set $V^o$ of strictly positive
dimension $m \leq n$. If $V^0$ does not intersect a closed
$n$-dimensional set $W$ of $V$ made up exclusively of all the
points of $V$ that are not in $V^o,$ then $V$ is reducible and
contains at least one irreducible component of dimension $m$.
\endproclaim

\demo {Proof} Without loss of generality, one can assume that the
$m$-dimensional variety $\overline {V^o}$ is irreducible\footnote
{Overline stands for the Zariski closure.}. Obviously $\overline
{V^o} \nsubseteq W$ since no point of an open set in $\overline
{V^o}$ lies in the closed set $W$. It must be so then that $V$
contains as a maximal irreducible component the $m$-dimensional
variety $\overline {V^o}$.
\enddemo
\enddemo

\head Section Two
\endhead

All the necessary preliminary material is now in place for a proof
of Theorem 1, and Theorem 5.

\proclaim {Theorem 1} Given an integer $k \ge 1 $ there exists an
infinite set $S_k$ consisting of parafree groups or rank 2 and
deviation $k$ having the property that for $G_i$, $G_j$ in $S_k$
the following holds: \roster
\item $Dim (R(G_i))=4+2k.$
\item If $i\ne j$ and $c=4+2k$, then $N_c(R(G_i))\ne N_c(R(G_j)).$
\item $R(G_i)\ncong R(G_j).$
\item Given any integer $m$ there exists an infinite subset $Q$ in
$S_k$ such that for any $G_t \in Q$ and $c$ as directly above
$N_c(R(G_t))\ge m$.
\endroster
\endproclaim

\demo {Proof} Without loss of generality $k\ge 1$ can be fixed.
Now, for each term, say the $j$-th term,  of the sequence
$kP_{i=1}^{\infty}$ the corresponding group $pGkP_j$ is a parafree
group or rank two by \cite {B1}, and since it maps onto a
non-trivial free product of cyclics involving $k+2$  factors it is
of deviation $k$.

\medskip
Now, refer to the Coarse Sieve Theorem of Section One, and set
$G_i=pGkP_i$, where on the right side is the group corresponding
to the $i$-term of the sequence $kP_i^{\infty}$. Set $\frak
N_i=N(x_{1}^2x_{2}^2)$ where the right side is the normal subgroup
of $pGkP_i$ generated by the closure of the word $x_{1}^2x_{2}^2$.
Using the outlined procedure, in the spirit of the Coarse Sieve
Theorem create an infinite sequence $S_{i=1}^{\infty}(G_i,\frak
N_i)$ of pairs $(G_1,\frak N_1), (G_2,\frak N_2), \cdots$
consisting of {\it fg} groups $G_i$ and a corresponding normal
subgroup $\frak N_i$ of $G_i$. Now let $N_c(V)$ be the {\it mirc}
counting function in dimension $c=4+2k$.

\medskip
By Lemma 1.5 it follows that $Dim(R(pGkP_i))=Dim(R(GkP_i))$ for
all $i\ge 1$. In other words, $Dim(R(G_i))=Dim(R(G_i/ \frak N_i))=
Dim(R(\frac {pGkP_i} {N(x_{1}^2x_{2}^2)}))=4+2k$ for all values of
$i$. Now, by Lemma 1.51 it follows that $N_c(R(\frac {pGkP_i}
{N(x_{1}^2x_{2}^2)}))$ assumes an infinite number of integral
points. In other words, $N_{4+2k}((R(G_i/ \frak N_i)))$ assumes an
infinite number of values. Thus the conditions demanded by the
Coarse Sieve Theorem are satisfied and all its consequences
follow. So employing the terminology of the Coarse Sieve Theorem
the sequence $S_{i=1}^{\infty}$ has an infinite subsequence $S^2$
such that for different pairs $(G_i,\frak N_i)$,\,\,$(G_j,\frak
N_j)$ in $S^2$ it is so that

$$N_c(R(G_i)) \ne N_c(R(G_j)). \tag {2-1}$$

\flushpar So the set $S_2=\{G \vert G\,\, \text { occurs in some
term of}\,\, S^2 \}$ can be ordered using the function $N_c(V)$.
Now for the sake of the statement of Theorem 1 rename $S_2$ by
$S_k$. Clearly, for $G_i$ and $G_j$ in $S_k$ it is the case that
$R(G_i)\ne R(G_j)$ since two algebraic varieties that are
isomorphic must have the same number of {\it mirc} of a given
dimension. Further, the infinite set $Q$ exists trivially since
$S_k$ has an infinite number of elements. This completes the proof
since $k$ was an arbitrary integer $\ge 1$.

\medskip The sequel will conclude with a statement and  proof of Theorem 5,
 a result instrumental in the proof of Theorem 3. Again, the prime
 sequence introduced in Section One, and the Coarse Sieve Theorem will prove
 indispensable.

 \proclaim {Theorem 5}
 Let $G_{pq}=<a,b,c; b^p[a,b^p]=c^q>$, where $p,q$ are different odd prime
 integer and where $[a,b]$ is defined to be
  $a^{-1}b^{-1}ab$. Then the following assertions hold.
 \roster
 \item $Dim(R(G_{pq}))=7.$

 \item $R(G_{pq})$ is reducible.

 \item There exists an infinite set $S$ of distinct prime
 integer pairs $S=\{(p,q)\vert p, q \in \Bbb Z_+ \}$ having the
 property that for any two different pairs $(p,q)$ and $(p',q')$ in $S$ the
 corresponding groups $G_{pq}$ and $G_{p'q'}$ have the property that
 $R(G_{pq})$ is not isomorphic  to
 $R(G_{p'q'})$. Also, the set $S$ can be well ordered using the {\it mirc}
 counting function $N_7(V)$ on representation varieties $V$ corresponding
 to groups $G_{pq}$ associated to $S$.

 \item The groups $G_{pq}$
  are residually torsion-free nilpotent.

 \item The groups $G_{pq}$ are parafree of rank 2 and deviation one.
 \endroster
 \endproclaim

 \demo {Proof}

 Proof 1: Let $<x>$ stand for the infinite cyclic group and
$\Bbb Z_t$ for the cyclic group of order $t$. Then the group
$<x>\ast \Bbb Z_p \ast \Bbb Z_q$ is a homomorphic image of the
group $G_{pq}$. Now $Dim(<x>\ast \Bbb Z_p \ast \Bbb Z_q)=7$; see
\cite {L1}. But, given any finitely generated group $G$ and  a
normal subgroup $N$ of $G$ it is easy to see that $R(G/N)\subseteq
R(G)$. So $Dim((R(G/N))\leq Dim(R(G))$. Consequently
$Dim(G_{pq})\ge 7$. That $Dim(G_{pq})\leq 7$ follows from the main
theorem of \cite {L1} which guarantees that any group with
presentation $G=<x_1,\dots,x_n, y;w=y^n>$, $n$ an integer,  and
$w$ a word not involving the generator $y$, has the property that
$Dim(R(G))\leq 3n+1$. Thus one is forced to conclude that
$Dim(R(G_{pq}))=7$.

\medskip
Proof 2: The reducibility of $R(G_{pq})$ stems from the fact that
if $H$ is a homomorphic image of $G$ of the same dimension as $G$
then $N_{Dim(R(G))}(R(G)) \ge N_{Dim(R(G))}(R(H))$. See \cite
{L4}. Also important is that the product of algebraic varieties is
reducible when at least one of the factors is a reducible variety.
But any finite non-trivial cyclic group has a reducible
representation variety \cite {L1}.

\medskip
Proof 3: The proof of part 3  employs the Coarse Sieve Theorem,
and the sequence $2P_{i=1}^\infty$ defined in the beginning of
Section One.
\medskip

\flushpar Associate with each term of $2P_{i=1}^\infty$  a pair
consisting of a group and a normal subgroup. For example with the
$k$-th term $2P_k$ of the sequence associate the pair

$$(G_{2P_k}, N_{2P_k})\tag {2-2}$$

\flushpar consisting of the group $G_{2P_k}=
G_{(P_{2k-1},P_{2k})}=<a,b,c;b^{P_{2k-1}}
[a,b^{P_{2k-1}}]=c^{P_{2k}}>$, and the normal subgroup $N_{2P_k}$
of $G_{2P_k}$ generated by the word  $\{b^{P_{2k-1}} \}$. Recall
that $P_{2k-1}, P_{2k}$ are successive primes in the sequence of
primes $P_{i=1}^\infty: 3,5,7,11,17,\dots$

\medskip\flushpar
It is straight forward to see that associated with any term of the
sequence $2P_{i=1}^\infty$,  say the $k$-th term  $2P_k$,  it is
so that for the corresponding pair $(G_{2P_k}, N_{2P_k})$ in (2-2)
the following holds:

$$G_{2P_k}/N_{2P_k}\cong <a,b,c; b^{P_{2k-1}}=1, c^{P_{2k}}=1>$$

\medskip
\flushpar Thus, since the representation variety of a free product
of a finite number of {\it
 fg} groups is the product of the corresponding
representation varieties of the factors, it follows that
$Dim(R(G_{2P_k}/N_{2P_k}))=7$. See \cite {L1}. So
$Dim(R(G_{2P_k}))=Dim R(G_{2P_k}/N_{2P_k})=7$.  Let $N_7(V)$ be
the {\it mirc} counting function in dimension 7, for any algebraic
variety $V$. \medskip

\flushpar Now, using elementary facts introduced in \cite {L1},
also in \cite {L3}, it follows immediately that
$N_7(R(G_{2P_k}/N_{2P_k})= \frac {(P_{2k-1}-1 )(P_{2k}-1)} {2^2}.$
Clearly, as $k$ tends to infinity $N_7(R(G_{2P_k}/N_{2P_k})$ also
tends to infinity. So the set $\{N_7(R(G_{2P_k}/N_{2P_k}))\vert
k\in 1,2,3,4,\dots \}$ has an infinite number of integer points.
The conditions of the CS Theorem are thus satisfied, and
consequently by wording the forth going in its language, one is
guaranteed an infinite set $S$ of distinct odd primes integers
pairs $S=\{(p,q)\vert p, q \in \Bbb Z_+ \}$ having the
 property that for different pair $(p,q)$ and $(p',q')$ in $S$ the
 corresponding groups $G_{pq}$ and $G_{p'q'}$ have the property that
 $R(G_{pq})$ is not isomorphic  to
 $R(G_{p'q'})$.

\medskip\flushpar
Proof 4:  In \cite {6},  G. Baumslag showed  that if
$G=<x_1,\dots,x_n, y; w(x_1,\dots,x_n)=y^t>$ is such that
$w(x_1,\cdots,x_n)\neq 1$ generates its own centralizer in
$<x_1,\cdots,x_n>$ and $t$ is a positive integer, that then $G$ is
residually torsion free nilpotent. Under the definition of $[a,b]$
stipulated in the statement of the theorem the word $b^p[a,b^p]$,
generates its own centralizer in the free group on $\{a,b\}$, and
consequently the groups $G_{pq}$ are residually torsion free
nilpotent.

\medskip\flushpar
Proof 5: G. Baumslag  in  \cite { B1} introduced a result
 quite handy in building non-isomorphic parafree groups  of the same
 rank as a previously given parafree group:
 \medskip\flushpar
 {\it
 Let $r$ and $n$ be positive integers and $H$ parafree of rank $r$, and
 let $(x)$ be the infinite cyclic group on $(x)$. Further, let
  $W\, \in \, H $. Suppose $W$ is the $k$ power of $W'$
 modulo $\gamma_2 H$, where $W'$ is itself not a power modulo
  $\gamma_2 H$; also, assume  that $k$ and $n$ are coprime and
 that $G$ is residually nilpotent. Then the  generalized free product
 $G = \{ H \ast (x) ; W = x^n \} $
 is  parafree  of rank $r$. }

 \medskip\flushpar Using Baumslag's result and Part 4, one can
 see that the groups $G_{pq}$ are parafree of rank 2. That they have
 deviation one stems from the fact that $Dim (R(G_{pq}))=7$ and
 thus the number of generators necessary to generate a $G_{pq}$
 is at least 3. So the deviation for each of the $G_{pq}$ is one.

 \enddemo

  \end